\documentclass[10pt,twoside,reqno]{amsart}
\usepackage{amssymb}
\usepackage{multirow}
\calclayout

\numberwithin{equation}{section}
\makeatletter

\renewcommand{\@secnumfont}{\bfseries}

\renewcommand{\section}{\@startsection{section}{1}%
  {0mm}{.7\linespacing\@plus\linespacing}{.5\linespacing}
  {\normalfont\bfseries\centering}}

\newcommand{\bibsection}{\@startsection{section}{1}%
  {0mm}{.7\linespacing\@plus\linespacing}{.5\linespacing}
  {\normalfont\scshape\centering}}

\renewcommand{\@biblabel}[1]{#1.}

\newtheorem{thm}{\bf Theorem}[section]

\newtheorem{cor}[thm]{\bf Corollary}

\begin{document}

\vspace{1.3cm}

\title {Some identities of degenerate ordered Bell polynomials and numbers arising from umbral calculus}

\author{Taekyun Kim}
\address{Department of Mathematics, College of Science Tianjin Polytechnic University, Tianjin 300160, China\\
Department of Mathematics, Kwangwoon University, Seoul 139-701, Republic
	of Korea}
\email{tkkim@kw.ac.kr}

\author{Dae San Kim}
\address{Department of Mathematics, Sogang University, Seoul 121-742, Republic
	of Korea}
\email{dskim@sogang.ac.kr}
\thanks{\scriptsize }

\maketitle
\markboth{\centerline{\scriptsize Degenerate ordered Bell polynomials and numbers associated with umbral calculus}}{\centerline{\scriptsize DAE SAN KIM AND TAEKYUN KIM}}

\begin{abstract}
In this paper, we study degenerate ordered Bell polynomials with the viewpoint of Carlitz's degenerate Bernoulli and Euler polynomials and derive by using umbral calculus some properties and new identities for the degenerate ordered Bell polynomials associated with special polynomials.

\end{abstract}

\bigskip
\medskip
\section{\bf Introduction}
As is well known, the ordinary Euler polynomials are defined by the generating function
\begin{equation}\begin{split}\label{01}
\frac{2}{e^t+1}e^{xt} = \sum_{n=0}^\infty E_n(x) \frac{t^n}{n},\quad (\textnormal{see} \,\, [6,7,8]).
\end{split}\end{equation}
When $x=0$, $E_n=E_n(0)$, $(n \geq 0)$, are called the Euler numbers. The Bernoulli polynomials are also given by the generating function as follows:
\begin{equation}\begin{split}\label{02}
\frac{t}{e^t-1}e^{xt} = \sum_{n=0}^\infty B_n(x) \frac{t^n}{n!},\quad (\textnormal{see} \,\, [5,6,8]).
\end{split}\end{equation}
When $x=0$, $B_n=B_n(0)$, are called the Bernoulli numbers. For $\lambda \in \mathbb{R}$, L. Carlitz considered the degenerate Euler and Bernoulli polynomials which are given by the generating function
\begin{equation}\begin{split}\label{03}
\frac{2}{(1+\lambda t)^{\frac{1}{\lambda}}+1}(1+\lambda t)^{\frac{x}{\lambda}} = \sum_{n=0}^\infty \mathcal{E}_n(x|\lambda) \frac{t^n}{n!},
\end{split}\end{equation}
and
\begin{equation}\begin{split}\label{04}
\frac{t}{(1+\lambda t)^{\frac{1}{\lambda}}-1}(1+\lambda t)^{\frac{x}{\lambda}} = \sum_{n=0}^\infty \beta_n(x|\lambda) \frac{t^n}{n!},\quad (\textnormal{see} \,\, [2]).
\end{split}\end{equation}
Note that $\lim_{\lambda \rightarrow 0} \mathcal{E}_n(x|\lambda) = E_n(x)$ and $\lim_{\lambda \rightarrow 0} \beta_n(x|\lambda) = B_n(x)$, $(n \geq 0)$, (see [2,7])
The falling factorial sequences are defined by
\begin{equation}\begin{split}\label{05}
(x)_0 =1,\,\,(x)_n = x(x-1)\cdots(x-(n-1)),\,\,(n \geq 1).
\end{split}\end{equation}
The Stirling numbers of the first kind are defined as
\begin{equation}\begin{split}\label{06}
(x)_n = \sum_{l=0}^n S_1(n,l) x^l,\,\,(n \geq 0),\quad (\textnormal{see} \,\, [10]).
\end{split}\end{equation}
The Stirling numbers of the second kind are also defined by
\begin{equation}\begin{split}\label{07}
x^n = \sum_{l=0}^n S_2(n,l) (x)_l,\,\,(n \geq 0),\quad (\textnormal{see} \,\, [6,10]).
\end{split}\end{equation}
It is well known that the ordered Bell polynomials are defined by the generating function
\begin{equation}\begin{split}\label{08}
\frac{1}{2-e^t}e^{xt} = \sum_{n=0}^\infty b_n(x) \frac{t^n}{n!},\quad (\textnormal{see} \,\, [3]).
\end{split}\end{equation}
When $x=0$, $b_n=b_n(0)$ are called the ordered Bell numbers. From \eqref{08}, we note that
\begin{equation}\begin{split}\label{09}
\frac{1}{2-e^t} &= \frac{1}{1-(e^t-1)} = \sum_{m=0}^\infty (e^t-1)^m \\
&= \sum_{m=0}^\infty m! \sum_{n=m}^\infty S_2(n,m) \frac{t^n}{n!} = \sum_{n=0}^\infty \left( \sum_{m=0}^n m! S_2(n,m) \right) \frac{t^n}{n!}.
\end{split}\end{equation}
Thus, by \eqref{09}, we get
\begin{equation}\begin{split}\label{10}
b_n = \sum_{m=0}^n m! S_2(n,m),\,\,(n \geq 0).
\end{split}\end{equation}
For $r \in \mathbb{N}$, the higher-order ordered Bell polynomials are given by the generating function
\begin{equation}\begin{split}\label{11}
\left( \frac{1}{2-e^t} \right)^r e^{xt} = \sum_{n=0}^\infty b_n^{(r)}(x) \frac{t^n}{n!}, \quad (\textnormal{see} \,\, [3]).
\end{split}\end{equation}
When $x=0$, $b_n^{(r)}=b_n^{(r)}(0)$, $(n \geq 0)$, are called the higher-order ordered Bell numbers.
By \eqref{11}, we easily get
\begin{equation}\begin{split}\label{12}
\sum_{n=0}^\infty b_n^{(r)} \frac{t^n}{n!} &= (2-e^t)^{-r} = \big(1-(e^t-1)\big)^{-r}\\
&= \sum_{m=0}^\infty {m+r-1 \choose m} (e^t-1)^m\\
&= \sum_{m=0}^\infty {m+r-1 \choose m} m! \sum_{n=m}^\infty S_2(n,m) \frac{t^n}{n!}\\
&= \sum_{n=0}^\infty \left( \sum_{m=0}^n {m+r-1 \choose m} m! S_2(n,m) \right) \frac{t^n}{n!}
\end{split}\end{equation}
From \eqref{12}, we note that
\begin{equation}\begin{split}\label{13}
b_n^{(r)} = \sum_{m=0}^n {m+r-1 \choose m} m! S_2(n,m),\,\,(n \geq 0).
\end{split}\end{equation}
In view of \eqref{03} and \eqref{04}, we consider the degenerate ordered Bell polynomials given by the generating function
\begin{equation}\begin{split}\label{14}
\frac{1}{2-(1+\lambda t)^{\frac{1}{\lambda }}}(1+\lambda t)^{\frac{x}{\lambda }} = \sum_{n=0}^\infty b_{n,\lambda} (x) \frac{t^n}{n!}.
\end{split}\end{equation}
By \eqref{14}, we easily get
\begin{equation}\begin{split}\label{15}
\lim_{\lambda  \rightarrow 0} b_{n,\lambda} (x) = b_n(x), \,\,(n \geq 0).
\end{split}\end{equation}
When $x=0$, $b_{n,\lambda }=b_{n,\lambda }(0)$ are called the degenerate ordered Bell numbers. Replacing $t$ by $\frac{1}{\lambda }\big(e^{\lambda t-1}\big)$, we get
\begin{equation}\begin{split}\label{16}
\sum_{m=0}^\infty b_{m,\lambda }(x) \frac{1}{m!} \left( \frac{1}{\lambda } \big(e^{\lambda t-1}\big) \right)^m &= \frac{1}{2-e^t}e^{xt}\\
&= \sum_{n=0}^\infty b_n(x) \frac{t^n}{n!}.
\end{split}\end{equation}
Now, we observe that
\begin{equation}\begin{split}\label{17}
&\sum_{m=0}^\infty b_{m,\lambda }(x) \frac{1}{m!} \lambda^{-m} \big( e^{\lambda t}-1\big)^m \\
&= \sum_{m=0}^\infty b_{m,\lambda }(x) \lambda^{-m} \sum_{n=m}^\infty S_2(n,m) \lambda^n \frac{t^n}{n!} \\
&= \sum_{n=0}^\infty \left( \sum_{m=0}^n \lambda^{n-m} b_{m,\lambda }(x) S_2(n,m) \right) \frac{t^n}{n!}.
\end{split}\end{equation}
From \eqref{16} and \eqref{17}, we have
\begin{equation}\begin{split}\label{18}
b_n(x) = \sum_{m=0}^n \lambda^{n-m} b_{m,\lambda }(x) S_2(n,m).
\end{split}\end{equation}

Here we study degenerate ordered Bell polynomials with the viewpoint of Carlitz's degenerate Bernoulli and Euler polynomials and derive by using umbral calculus some properties and new identities for the degenerate ordered Bell polynomials associated with special polynomials.

\section{\bf Quick review of umbral calculus}
Let $\mathbb{C}$ be the complex number field and let $\mathcal{F}$ be the set of all formal power series with coefficients in $\mathbb{C}$ in the variable $t$:
\begin{equation}\begin{split}\label{19}
\mathcal{F}= \left\{ f(t)= \sum_{k=0}^\infty a_k \frac{t^k}{k!}\,\,   \bigg| \,\,a_k \in \mathbb{C}    \right\}.
\end{split}\end{equation}
We denote the algebra $\mathbb{C}[x]$ of polynomials in $x$ over the field $\mathbb{C}$ by $\mathbb{P}$ ($\mathbb{P} = \mathbb{C}[x]$). Let $\mathbb{P}^*$ be the vector space of all linear functionals on $\mathbb{P}$, and let $<L \,|\,p(x)>$ denote the action of the linear functional $L$ on $p(x)$, which satisfies $<L+M|p(x)> = <L|p(x)> + <M|p(x)>$ and $<cL|p(x)> = c<L|p(x)>$, where $c$ is a complex number. The linear functional $<f(t)| \cdot>$ on $\mathbb{P}$ is defined by $<f(t)|x^n>=a_n$ $(n \geq 0)$, where $f(t)= \sum_{k=0}^\infty a_k \frac{t^k}{k!} \in \mathcal{F}$. Thus we note that $<t^k|x^n> = n! \delta_{n,k}$ $(n,k \geq 0)$, where $\delta_{n,k}$ is the Kronecker symbol (see [1,4,10]). Let $f_L(t) = \sum_{k=0}^\infty \frac{<L|x^k>}{k!}t^k$. Then we have $<f_L(t)|x^n> = <L|x^n>$ $(n \geq 0)$. So, the map $L \longmapsto f_L(t)$ is a vector space isomorphism from $\mathbb{P}^*$ onto $\mathcal{F}$. Henceforth, $\mathcal{F}$ denotes both the algebra of formal power series in $t$ and the vector space of all linear functionals on $\mathbb{P}$, and so an element $f(t)$ of $\mathcal{F}$ will be thought of as both a formal power series and a linear functional. We call $\mathcal{F}$ the umbral algebra.  The umbral calculus is the study of the umbral algebra. Let $f(t) (\neq 0) \in \mathcal{F}$. Then order of $f(t)$ is the smallest positive integer $k$ for which the coefficient of $t^k$ does not vanish. The order of $f(t)$ is denotes by $o(f(t))$, (see [7,8,10]). For $f(t) \in \mathcal{F}$ and $p(x) \in \mathbb{P}$, we have
\begin{equation}\begin{split}\label{20}
f(t) = \sum_{k=0}^\infty <f(t) | x^k > \frac{t^k}{k!}, \,\, p(x) = \sum_{k=0}^\infty <t^k | p(x)>  \frac{t^k}{k!}.
\end{split}\end{equation}
Thus, by \eqref{20}, we easily get
\begin{equation}\begin{split}\label{21}
p^{(k)}(0) = <t^k | p(x)> = <1|p^{(k)}(x)>, \quad (\textnormal{see} \,\, [9,10]),
\end{split}\end{equation}
where $p^{(k)}(x) = \left( \frac{d}{dx} \right)^k p(x)$. From \eqref{21}, we note that
\begin{equation}\begin{split}\label{22}
t^k p(x) = p^{(k)}(x), \,\, e^{yt}p(x) = p(x+y), \,\,\text{and}\,\, <e^{yt}|p(x)>=p(y),\quad (\textnormal{see} \,\, [10,11,12]).
\end{split}\end{equation}
For $f(t),g(t) \in \mathcal{F}$ with $o(f(t))=1$, $o(g(t))=0$, there exists a unique sequence $S_n(x)$ of polynomials with $\deg S_n(x)=n$ such that
\begin{equation}\begin{split}\label{23}
<g(t)f(t)^k | S_n(x) > = n! \delta_{n,k},\,\,(n,k \geq 0),\quad (\textnormal{see} \,\, [6,7,10]).
\end{split}\end{equation}
The sequence $S_n(x)$ is called the Sheffer sequence for $\big(g(t),f(t)\big)$, and we write $S_n(x)\sim \big( g(t), f(t)\big)$. It is well known that
\begin{equation*}\begin{split}
S_n(x)\sim \big( g(t), f(t)\big) \Longleftrightarrow \frac{1}{g\big(\bar{f}(t)\big)}e^{x\bar{f}(t)} = \sum_{n=0}^\infty S_n(x) \frac{t^n}{n!},
\end{split}\end{equation*}
where $\bar{f}(t)$ is the compositional inverse of $f(t)$ such that $f\big( \bar{f}(t)\big) = \bar{f} \big( f(t) \big) = t$. Let $S_n(x)\sim \big( g(t), f(t)\big)$. Then we have
\begin{equation}\begin{split}\label{24}
f(t)s_n(x) = nS_{n-1}(x),\,\,(n\geq 1), \,\, S_n(x) = \sum_{j=0}^n \frac{<g\big(\bar{f}(t)\big)^{-1}\bar{f}(t)^j|x^n>}{j!} x^j,
\end{split}\end{equation}
\begin{equation}\begin{split}\label{25}
S_n(x+y) = \sum_{j=0}^n {n \choose j} S_j(x) P_{n-j}(y), \,\,\text{where}\,\,P_n(y) = g(t)S_n(y),
\end{split}\end{equation}
\begin{equation}\begin{split}\label{26}
<f(t)|xp(x)> = <\delta_t f(t)|p(x)>,\,\, \text{where}\,\, \delta_t f(t) = \frac{d}{dt}f(t),
\end{split}\end{equation}
and
\begin{equation}\begin{split}\label{27}
S_{n+1}(x) = \left( x- \frac{g'(t)}{g(t)} \right) \frac{1}{f'(t)}S_n(x),\,\,(n \geq 0), \quad (\textnormal{see} \,\, [10]).
\end{split}\end{equation}
For $p_n(x) \sim \big( 1, f(t) \big)$, $q_n(x) \sim \big(1, g(t) \big)$, we have
\begin{equation}\begin{split}\label{28}
q_n(x) = x \left( \frac{f(t)}{g(t)} \right)^n x^{-1} p_n(x),\,\,(n \geq 1), \quad (\textnormal{see} \,\, [10,12]).
\end{split}\end{equation}
Let us consider the following two Sheffer sequences:
\begin{equation}\begin{split}\label{29}
S_n(x) \sim \big(g(t),f(t) \big),\,\, r_n(x) \sim \big( h(t), l(t) \big).
\end{split}\end{equation}
Then we have
\begin{equation}\begin{split}\label{30}
S_n(x) = \sum_{m=0}^n C_{n,m} r_m(x), \,\,(n \geq 0),
\end{split}\end{equation}
where
\begin{equation}\begin{split}\label{31}
C_{n,m} = \frac{1}{m!} \Big< \frac{h \big( \bar{f}(t)\big)}{g\big(\bar{f}(t)\big)}\Big(l(\big(\bar{f}(t)\big)\Big)^m | x^n \Big>,\quad (\textnormal{see} \,\, [10]).
\end{split}\end{equation}
Finally, we recall that
\begin{equation}\label{100}
\frac{d}{dx}S_n(x)=\sum_{k=0}^{n-1} \binom{n}{k}<\bar{f}(t)|x^{n-k}>S_k(x).
\end{equation}

\section{Degenerate ordered Bell numbers and polynomials associated with umbral calculus}

From \eqref{14}, we note that the degenerate ordered Bell polynomials are defined by the generating function
\begin{equation*}\begin{split}
\frac{1}{2-(1+\lambda t)^{\frac{1}{\lambda }}}(1+\lambda t)^{\frac{x}{\lambda }} = \sum_{n=0}^\infty b_{n,\lambda }(x) \frac{t^n}{n!}.
\end{split}\end{equation*}
When $x=0$, $b_{n,\lambda }=b_{n,\lambda }(0)$ are called the ordered Bell numbers. By \eqref{14}, we get
\begin{equation}\begin{split}\label{32}
&\sum_{n=0}^\infty b_{n,\lambda } \frac{t^n}{n!}=\frac{1}{2-(1+\lambda t)^{\frac{1}{\lambda }}}= \frac{1}{1-\big((1+\lambda t)^{\frac{1}{\lambda }}-1\big)}\\
&= \sum_{m=0}^\infty \left( (1+\lambda t)^{\frac{1}{\lambda }}-1 \right)^m = \sum_{m=0}^\infty \left( e^{\frac{1}{\lambda }\log(1+\lambda t)}-1 \right)^m\\
&=\sum_{m=0}^\infty m! \sum_{k=m}^\infty S_2(k,m) \lambda^{-k} \frac{1}{k!} \big(\log(1+\lambda t)\big)^k \\
&= \sum_{k=0}^\infty \left( \sum_{m=0}^k m! S_2(k,m) \lambda^{-k} \right) \sum_{n=k}^\infty S_1(n,k) \frac{\lambda ^n t^n}{n!}\\
&= \sum_{n=0}^\infty \left( \sum_{k=0}^n \sum_{m=0}^k m! S_2(k,m) S_1(n,k) \lambda^{n-k} \right) \frac{t^n}{n!}.
\end{split}\end{equation}
From \eqref{32}, we note that
\begin{equation}\begin{split}\label{33}
b_{n,\lambda } =  \sum_{k=0}^n \sum_{m=0}^k m! S_2(k,m) S_1(n,k) \lambda^{n-k}, \,\,(n \geq 0).
\end{split}\end{equation}
By \eqref{14} and \eqref{23}, we get
\begin{equation}\begin{split}\label{34}
b_{n,\lambda }(x) \sim \big( 2-e^t, \tfrac{1}{\lambda }(e^{\lambda t}-1)\big),\,\,(n \geq 0).
\end{split}\end{equation}
That is
\begin{equation}\begin{split}\label{35}
\frac{1}{2-e^{\frac{1}{\lambda }\log(1+\lambda t)}} e^{x \frac{1}{\lambda } \log(1+\lambda t)} = \sum_{ n=0}^\infty b_{n,\lambda }(x) \frac{t^n}{n!}.
\end{split}\end{equation}
Let $f(t)$ be the linear functional such that
\begin{equation}\begin{split}\label{36}
<f(t)|p(x)> = \int_0^y p(u)du.\,\,\text{for all}\,\,p(x).
\end{split}\end{equation}
Then, by \eqref{20} and \eqref{36}, we get
\begin{equation}\begin{split}\label{37}
f(t) &= \sum_{k=0}^\infty \frac{<f(t)|x^k>}{k!}t^k = \sum_{k=0}^\infty \frac{y^{k+1}}{(k+1)!}t^k\\
&= \frac{1}{t}(e^{yt}-1).
\end{split}\end{equation}
From \eqref{36}, we have
\begin{equation}\begin{split}\label{38}
<\tfrac{1}{t}(e^{yt}-1)|p(x)> = \int_0^y p(u)du,
\end{split}\end{equation}
and
\begin{equation}\begin{split}\label{39}
\tfrac{1}{t}(e^{yt}-1) p(x) = \int_x^{x+y}p(u)du.
\end{split}\end{equation}
By \eqref{32}, \eqref{35}, we get
\begin{equation}\begin{split}\label{40}
&\sum_{n=0}^\infty b_{n,\lambda }(x) \frac{t^n}{n!} = \frac{1}{2-(1+\lambda t)^{\frac{1}{\lambda }}} e^{\frac{x}{\lambda }\log(1+\lambda t)}\\
&= \left( \sum_{l=0}^\infty b_{l,\lambda }\frac{t^l}{l!} \right) \left( \sum_{m=0}^\infty \left(\frac{x}{\lambda}\right)^m \frac{1}{m!} \big(\log(1+\lambda t)\big)^m\right)\\
&=\left( \sum_{l=0}^\infty b_{l,\lambda }\frac{t^l}{l!} \right) \left( \sum_{m=0}^\infty \left(\frac{x}{\lambda}\right)^m \sum_{k=m}^\infty S_1(k,m) \frac{\lambda ^k t^k}{k!} \right)\\
&=\left( \sum_{l=0}^\infty b_{l,\lambda }\frac{t^l}{l!} \right) \left( \sum_{k=0}^\infty \left( \sum_{m=0}^k \lambda ^{k-m} S_1(k,m) x^m \right) \frac{t^k}{k!} \right)\\
&=\sum_{n=0}^\infty \left\{ \sum_{k=0}^n \sum_{m=0}^k \lambda ^{k-m} S_1(k,m) {n \choose k} x^m b_{n-k,\lambda } \right\} \frac{t^n}{n!}.
\end{split}\end{equation}
Comparing the coefficients on both sides of \eqref{40}, we have
\begin{equation}\begin{split}\label{41}
b_{n,\lambda }(x) = \sum_{k=0}^n \sum_{m=0}^k {n \choose k}\lambda ^{k-m} S_1(k,m)b_{n-k,\lambda } x^m .
\end{split}\end{equation}
From \eqref{39}, we note that
\begin{equation}\begin{split}\label{42}
&\frac{e^t-1}{t}b_{n,\lambda }(x) = \int_x^{x+1} b_{n,\lambda }(u)du\\
&=\sum_{k=0}^n  \sum_{m=0}^k {n \choose k} \lambda ^{k-m} S_1(k,m) b_{n-k,\lambda } \int_x^{x+1}u^m du\\
&=\sum_{k=0}^n  \sum_{m=0}^k {n \choose k} \lambda ^{k-m} S_1(k,m) b_{n-k,\lambda } \frac{1}{m+1}\big( (x+1)^{m+1}-x^{m+1}\big)\\
&= \sum_{k=0}^n \sum_{m=0}^k \sum_{l=0}^m {m+1 \choose l} {n \choose k} \lambda ^{k-m} S_1(k,m) b_{n-k,\lambda } \frac{1}{m+1}x^l.
\end{split}\end{equation}
As is well knwon, $B_n(x) = \frac{t}{e^t-1}x^n$, $(n \geq 0)$. By \eqref{42}, we get
\begin{equation}\begin{split}\label{43}
b_{n,\lambda }(x) = \sum_{k=0}^n \sum_{m=0}^k \sum_{l=0}^m  \lambda ^{k-m}{m+1 \choose l} {n \choose k}S_1(k,m)b_{n-k,\lambda } \frac{1}{m+1} B_l(x).
\end{split}\end{equation}
From \eqref{43}, we note that
\begin{equation}\begin{split}\label{44}
tb_{n,\lambda }(x) &= \sum_{k=1}^n \sum_{m=1}^k \sum_{l=1}^m  \lambda ^{k-m}{m+1 \choose l} {n \choose k} \frac{l}{m+1}  B_{l-1}(x) S_1(k,m)b_{n-k,\lambda}\\
&= \sum_{k=1}^n \sum_{m=1}^k \sum_{l=1}^m  \lambda ^{k-m}{m \choose l-1} {n \choose k} B_{l-1}(x) S_1(k,m)b_{n-k,\lambda}.
\end{split}\end{equation}
and
\begin{equation}\begin{split}\label{45}
&\frac{e^{yt}-1}{t}b_{n,\lambda }(x) = \int_x^{x+y} b_{n,\lambda }(u)du\\
&=\sum_{k=0}^n \sum_{m=0}^k \sum_{l=0}^m\lambda ^{k-m}{m+1 \choose l} {n \choose k}S_1(k,m)b_{n-k,\lambda } \left( \frac{B_{l+1}(x+y)-B_{l+1}(x)}{(m+1)(l+1)} \right)\\
&=\sum_{k=0}^n \sum_{m=0}^k \sum_{l=0}^m \sum_{j=0}^{l}  {m+1 \choose l} {n \choose k} {l+1 \choose j} \lambda ^{k-m}S_1(k,m)b_{n-k,\lambda }y^{l+1-j}B_j(x) \frac{1}{(m+1)(l+1)}\\
&=\sum_{k=0}^n \sum_{m=0}^k \sum_{l=0}^m \sum_{j=0}^l \frac{{m+1 \choose l} {n \choose k} {l+1 \choose j} }{(m+1)(l+1)  } \lambda ^{k-m}S_1(k,m)b_{n-k,\lambda }y^{l+1-j}B_j(x).
\end{split}\end{equation}
By multiplying $t$ on both sides of \eqref{45}, we get
\begin{equation}\begin{split}\label{46}
&b_{n,\lambda }(x+y)-b_{n,\lambda }(x) = \big(e^{yt}-1\big) b_{n,\lambda }(x) \\
&=\sum_{k=1}^n \sum_{m=1}^k \sum_{l=1}^m \sum_{j=1}^l \frac{{m+1 \choose l} {n \choose k} {l+1 \choose j}j }{(m+1)(l+1)  }j\lambda ^{k-m}S_1(k,m)b_{n-k,\lambda }y^{l+1-j}B_{j-1}(x)\\
&=\sum_{k=1}^n \sum_{m=1}^k \sum_{l=1}^m \sum_{j=1}^l \frac{{m+1 \choose l} {n \choose k} {l \choose j-1} }{m+1} \lambda ^{k-m}S_1(k,m)b_{n-k,\lambda }y^{l+1-j}B_{j-1}(x).
\end{split}\end{equation}
Therefore, by \eqref{46}, we obtain the following theorem.
\begin{thm}
For $n \geq 1$, we have
\begin{equation*}\begin{split}
&b_{n,\lambda }(x+y)-b_{n,\lambda }(x) \\
&=\sum_{k=1}^n \sum_{m=1}^k \sum_{l=1}^m \sum_{j=1}^l \frac{{m+1 \choose l} {n \choose k} {l \choose j-1} }{m+1} \lambda ^{k-m}S_1(k,m)b_{n-k,\lambda }y^{l+1-j}B_{j-1}(x).
\end{split}\end{equation*}
\end{thm}
Now, we observe that
\begin{equation}\begin{split}\label{47}
&\sum_{n=0}^\infty \big( -b_{n,\lambda }(x+1)+2b_{n,\lambda }(x) \big) \frac{t^n}{n!} = \frac{(1+\lambda t)^{\frac{x}{\lambda }}\big(2-(1+\lambda t)^{\frac{1}{\lambda }}\big)}{2-(1+\lambda t)^{\frac{1}{\lambda }}}\\
&= (1+\lambda t)^{\frac{x}{\lambda}} = \sum_{k=0}^\infty \left( \frac{x}{\lambda } \right)^k \frac{1}{k!} \big( \log(1+\lambda t) \big)^k\\
&= \sum_{k=0}^\infty \left( \frac{x}{\lambda } \right)^k \sum_{n=k}^\infty S_1(n,k) \lambda ^n \frac{t^n}{n!} = \sum_{n=0}^\infty \left( \sum_{k=0}^n \lambda ^{n-k} S_1(n,k) x^k \right) \frac{t^n}{n!}.
\end{split}\end{equation}
Comparing the coefficients on both sides of \eqref{47}, we obtain the following theorem.
\begin{thm}
For $n \geq 0$, we have
\begin{equation*}\begin{split}
2b_{n,\lambda }(x)-b_{n,\lambda }(x+1) = \sum_{k=0}^n \lambda ^{n-k} S_1(n,k) x^k.
\end{split}\end{equation*}
\end{thm}
From Theorem 1 and Theorem 2, we note that
\begin{equation}\begin{split}\label{48}
&b_{n,\lambda }(x+1) -b_{n,\lambda }(x) = b_{n,\lambda }(x) - \sum_{k=0}^n \lambda ^{n-k} S_1(n,k) x^k\\
&= \sum_{k=1}^n \sum_{m=1}^k \sum_{l=1}^m \sum_{j=1}^l \frac{{m+1 \choose l}{n \choose k}{l \choose j-1}}{m+1} \lambda^{k-m} S_1(k,m) b_{n-k,\lambda } B_{j-1}(x).
\end{split}\end{equation}
Therefore, by \eqref{48}, we obtain the following corollary.
\begin{cor}
For $n \geq 1$, we have
\begin{equation*}\begin{split}
& b_{n,\lambda }(x) = \sum_{k=1}^n \lambda ^{n-k}S_1(n,k)x^k\\
&+ \sum_{k=1}^n \sum_{m=1}^k \sum_{l=1}^m \sum_{j=1}^l \frac{{m+1 \choose l}{n \choose k}{l \choose j-1}}{m+1} \lambda^{k-m} S_1(k,m) b_{n-k,\lambda } B_{j-1}(x).
\end{split}\end{equation*}
\end{cor}
From $<t^k|x^n>=n! \delta_{n,k}$, $(n,k \geq 0)$, we note that
\begin{equation}\begin{split}\label{49}
b_{n,\lambda }(y)&= \Big< \sum_{l=0}^\infty b_{l,\lambda }(y) \frac{t^l}{l!} \Big| x^n \Big> = \Big< \frac{1}{2-(1+\lambda t)^{\frac{1}{\lambda }}}(1+\lambda t)^{\frac{y}{\lambda}}\Big|x^n \Big>\\
&= \Big< \frac{1}{2-(1+\lambda t)^{\frac{1}{\lambda }}} \Big| (1+\lambda t)^{\frac{y}{\lambda }}x^n \Big> \\
&= \sum_{m=0}^n {n \choose m} \sum_{k=0}^m y^k \lambda ^{m-k}S_1(m,k) \Big< \frac{1}{2-(1+\lambda t)^{\frac{1}{\lambda }}} \Big| x^{n-m} \Big> \\
&=  \sum_{m=0}^n {n \choose m} \sum_{k=0}^m y^k \lambda ^{m-k}S_1(m,k) \sum_{l=0}^{n-m} \frac{b_{l,\lambda }}{l!} <t^l | x^{n-m}>\\
&= \sum_{m=0}^n \sum_{k=0}^m {n \choose m} \lambda ^{m-k} S_1(m,k) y^k b_{n-m,\lambda }.
\end{split}\end{equation}
Therefore, by \eqref{49}, we obtain the following theorem.
\begin{thm}
For $n \geq 0$, we have
\begin{equation*}\begin{split}
b_{n,\lambda }(x) =  \sum_{m=0}^n \sum_{k=0}^m {n \choose m} \lambda ^{m-k} S_1(m,k) b_{n-m,\lambda}x^k.
\end{split}\end{equation*}
\end{thm}
We easily see that
\begin{equation}\begin{split}\label{50}
x^n \sim (1,t), \,\, (2-e^t)b_{n,\lambda }(x) \sim \big(1, \tfrac{1}{\lambda } (e^{\lambda t}-1) \big).
\end{split}\end{equation}
Then, by \eqref{28} and \eqref{50}, we get
\begin{equation}\begin{split}\label{51}
(2-e^t)b_{n,\lambda }(x) &= x \left( \frac{\lambda t}{e^{\lambda t}-1} \right)^n x^{-1} x^n\\
&= x \sum_{l=0}^\infty \lambda ^l \frac{B_l^{(n)}}{l!} t^l x^{n-1}\\
&= x \sum_{l=0}^{n-1} \lambda ^l {n-1 \choose l} B_l^{(n)} x^{n-1-l} = \sum_{l=0}^{n-1} \lambda ^l {n-1 \choose l} B_l^{(n)} x^{n-l},
\end{split}\end{equation}
where $B_n^{(\alpha)}$ are the higher-order Bernoulli numbers defined by the generating function
\begin{equation*}\begin{split}
\left( \frac{t}{e^t-1} \right)^\alpha = \sum_{n=0}^\infty B_n^{(\alpha)} \frac{t^n}{n!}.
\end{split}\end{equation*}
Thus, by \eqref{51}, we get
\begin{equation}\begin{split}\label{52}
b_{n,\lambda }(x) &= \sum_{l=0}^{n-1} \lambda ^l {n-1 \choose l} B_l^{(n)} \frac{1}{2-e^t}x^{n-l}\\
&= \sum_{l=0}^{n-1} \lambda ^l {n-1 \choose l} B_l^{(n)} b_{n-l}(x),\,\,(n \geq 1).
\end{split}\end{equation}
Therefore, by \eqref{52}, we obtain the following theorem.
\begin{thm}
For $n \geq 1$, we have
\begin{equation*}\begin{split}
b_{n,\lambda}(x) &= \sum_{l=0}^{n-1} \lambda ^l {n-1 \choose l} B_l^{(n)} b_{n-l,\lambda }(x)\\
&= \sum_{l=0}^n \sum_{m=0}^{n-l} \lambda ^l {n-1 \choose l} {n-l \choose m} B_l^{(n)} b_m x^{n-l-m}.
\end{split}\end{equation*}
\end{thm}
Let $S_n(x) \sim \big( 1, \tfrac{1}{\lambda }(e^{\lambda t}-1) \big)$. Then, by \eqref{23}, we get
\begin{equation}\begin{split}\label{53}
\sum_{n=0}^\infty S_n(x) \frac{t^n}{n!} &= e^{\frac{x}{\lambda }\log(1+\lambda t)} = (1+\lambda t)^{\frac{x}{\lambda }}\\
&= \sum_{n=0}^\infty \left( \frac{x}{\lambda } \right)_n \lambda ^n \frac{t^n}{n!}= \sum_{n=0}^\infty \frac{x}{\lambda }\left(\frac{x}{\lambda }-1 \right) \cdots \left( \frac{x}{\lambda }-n+1 \right) \lambda ^n \frac{t^n}{n!}\\
&=\sum_{n=0}^\infty x(x-\lambda )(x-2\lambda )\cdots(x-(n-1)\lambda ) \frac{t^n}{n!}\\
&= \sum_{n=0}^\infty (x)_{n,\lambda } \frac{t^n}{n!},
\end{split}\end{equation}
where $(x)_{0,\lambda }=1$, $(x)_{n,\lambda }=x(x-\lambda )(x-2\lambda )\cdots (x-(n-1)\lambda )$, $(n \geq 1)$. From \eqref{53}, we note that $(x)_{n,\lambda } \sim \big(1,\tfrac{1}{\lambda }(e^{\lambda t}-1)\big)$. Now, we observe that
\begin{equation}\begin{split}\label{54}
\sum_{n=0}^\infty (x)_{n,\lambda } \frac{t^n}{n!}&= (1+\lambda t)^{\frac{x}{\lambda }} = e^{\frac{x}{\lambda }\log(1+\lambda t)}\\
&= \sum_{m=0}^\infty \left( \frac{x}{\lambda }\right)^m \frac{1}{m!} \big( \log(1+\lambda t)\big)^m\\
&= \sum_{n=0}^\infty \left( \sum_{m=0}^n \lambda ^{n-m} S_1(n,m) x^m \right) \frac{t^n}{n!}.
\end{split}\end{equation}
Comparing the coefficients on both sides of \eqref{54}, we have
\begin{equation}\begin{split}\label{55}
(x)_{n,\lambda } = \sum_{m=0}^n \lambda ^{n-m} S_1(n,m) x^m,\,\,(n \geq 0).
\end{split}\end{equation}
For $(2-e^t)b_{n,\lambda }(x) \sim \big(1, \tfrac{1}{\lambda }(e^{\lambda t}-1)\big)$, $(x)_{n,\lambda } \sim \big(1, \tfrac{1}{\lambda }(e^{\lambda t}-1)\big)$, we get
\begin{equation}\begin{split}\label{56}
(2-e^t)b_{n,\lambda }(x) = (x)_{n,\lambda }  = \sum_{m=0}^n \lambda ^{n-m} S_1(n,m) x^m.
\end{split}\end{equation}
By \eqref{56}, we easily see that
\begin{equation}\begin{split}\label{57}
b_{n,\lambda }(x) &= \sum_{m=0}^n \lambda ^{n-m} S_1(n,m) \frac{1}{2-e^t} x^m\\
&=\sum_{m=0}^n \lambda ^{n-m} S_1(n,m) b_m(x),\,\,(n \geq 0).
\end{split}\end{equation}
\begin{cor}
For $n \geq 0$, we have
\begin{equation*}\begin{split}
b_{n,\lambda }(x) = \sum_{m=0}^n \lambda ^{n-m} S_1(n,m) b_m(x).
\end{split}\end{equation*}
\end{cor}
By \eqref{27} and \eqref{34}, we get
\begin{equation}\begin{split}\label{58}
b_{n+1,\lambda }(x) &= \left( x - \frac{(-e^t)}{2-e^t} \right) \frac{1}{e^{\lambda t}} b_{n,\lambda }(x) \\
&=\left( x+ \frac{e^t}{2-e^t} \right) b_{n,\lambda }(x-\lambda ) = xb_{n,\lambda }(x-\lambda ) + \frac{e^t}{2-e^t} b_{n,\lambda }(x-\lambda )\\
&= xb_{n,\lambda }(x-\lambda ) + \frac{e^t-2+2}{2-e^t}b_{n,\lambda }(x-\lambda )\\
&= xb_{n,\lambda }(x-\lambda )-b_{n,\lambda }(x-\lambda )+\frac{2}{2-e^t}b_{n,\lambda }(x-\lambda ).
\end{split}\end{equation}
From \eqref{41}, we note that
\begin{equation}\begin{split}\label{59}
& \frac{1}{2-e^t}b_{n,\lambda }(x-\lambda )= \frac{1}{2-e^t} \sum_{k=0}^n \sum_{m=0}^k {n \choose k}\lambda ^{k-m} S_1(k,m) b_{n-k,\lambda }(x-\lambda )^m\\
&= \frac{1}{2-e^t} \sum_{k=0}^n \sum_{m=0}^k {n \choose k} \lambda ^{k-m} S_1(k,m) b_{n-k,\lambda } \sum_{l=0}^m {m \choose l} (-1)^l \lambda ^l x^{m-l}\\
&= \sum_{k=0}^n \sum_{m=0}^k \sum_{l=0}^m {n \choose k} {m \choose l} \lambda ^{k-m+l} S_1(k,m) b_{n-k,\lambda }(-1)^l \frac{1}{2-e^t}x^{m-l}\\
&=  \sum_{k=0}^n \sum_{m=0}^k \sum_{l=0}^m {n \choose k} {m \choose l} \lambda ^{k-m+l} S_1(k,m) b_{n-k,\lambda }(-1)^l b_{m-l}(x).
\end{split}\end{equation}
By \eqref{58} and \eqref{59}, we get
\begin{equation}\begin{split}\label{60}
&b_{n+1,\lambda }(x) = xb_{n,\lambda }(x-\lambda )-b_{n,\lambda }(x-\lambda )+\frac{2}{2-e^t}b_{n,\lambda }(x-\lambda )\\
&= xb_{n,\lambda }(x-\lambda )-b_{n,\lambda }(x-\lambda )\\
&\quad +2\sum_{k=0}^n \sum_{m=0}^k \sum_{l=0}^m {n \choose k} {m \choose l} \lambda ^{k-m+l} S_1(k,m) b_{n-k,\lambda }(-1)^l b_{m-l}(x).
\end{split}\end{equation}
Therefore, by \eqref{60}, we obtain the following theorem.

\begin{thm}
For $n \geq 0$, we have
\begin{equation*}\begin{split}
&b_{n+1,\lambda }(x) - xb_{n,\lambda }(x-\lambda )+b_{n,\lambda }(x-\lambda )\\
&=2\sum_{k=0}^n \sum_{m=0}^k \sum_{l=0}^m {n \choose k} {m \choose l} \lambda ^{k-m+l} S_1(k,m) b_{n-k,\lambda }(-1)^l b_{m-l}(x).
\end{split}\end{equation*}
\end{thm}

\noindent{\bf{Remark.}} From $b_{n,\lambda }(x) \sim \big(2-e^t, \tfrac{1}{\lambda }(e^{\lambda t}-1)\big)$ and \eqref{100}, we note that
\begin{equation}\begin{split}\label{61}
\frac{d}{dx}b_{n,\lambda }(x) &= \sum_{l=0}^{n-1} {n \choose l} <\frac{1}{\lambda}\log (1+\lambda t)|x^{n-l}> b_{l,\lambda }(x)\\
 &= \sum_{l=0}^{n-1} {n \choose l} <\sum_{m=1}^{\infty}\frac{(-\lambda)^{m-1} t^m}{m}|x^{n-l}> b_{l,\lambda }(x)\\
&=\sum_{l=0}^{n-1}{n \choose l}(-\lambda)^{n-l-1}(n-l-1)!  b_{l,\lambda }(x), \,\,(n \geq 1).
\end{split}\end{equation}
Let
\begin{equation}\begin{split}\label{62}
\mathbb{P}_n = \left\{ p(x) \in \mathbb{C}[x] | \deg p(x) \leq n \right\},\,\,(n \geq 0).
\end{split}\end{equation}
For $p(x) \in \mathbb{P}_n$, let
\begin{equation}\begin{split}\label{63}
p(x) = \sum_{l=0}^n a_l b_{l,\lambda }(x).
\end{split}\end{equation}
From \eqref{34}, we note that
\begin{equation}\begin{split}\label{64}
<(2-e^t)(\tfrac{1}{\lambda }(e^{\lambda t}-1))^m|b_{n,\lambda }(x)> = n! \delta_{n,m}, \,\,(n, m\geq 0).
\end{split}\end{equation}
By \eqref{63} and \eqref{64}, we get, for $0\leq m\leq n$,
\begin{equation}\begin{split}\label{65}
&<(2-e^t)(\tfrac{1}{\lambda }(e^{\lambda t}-1))^m|p(x)>\\
&=\sum_{l=0}^n a_l < (2-e^t) (\tfrac{1}{\lambda }(e^{\lambda t}-1))^m | b_{l,\lambda }(x)>\\
&=\sum_{l=0}^n a_l l! \delta_{l,m} = a_m m!.
\end{split}\end{equation}
From \eqref{65}, we have
\begin{equation}\begin{split}\label{66}
a_m = \frac{1}{m!} <(2-e^t)(\tfrac{1}{\lambda }(e^{\lambda t}-1))^m|p(x)>,
\end{split}\end{equation}
where $n \geq m \geq 0$. Therefore, we obtain the following theorem.
\begin{thm}
For $p(x) \in \mathbb{P}_n$, we have
\begin{equation*}\begin{split}
p(x) = \sum_{m=0}^n a_m b_{m,\lambda }(x),\,\,(n \geq 0),
\end{split}\end{equation*}
where
\begin{equation*}\begin{split}
a_m = \frac{1}{m!} <(2-e^t)(\tfrac{1}{\lambda }(e^{\lambda t}-1))^m|p(x)>.
\end{split}\end{equation*}
\end{thm}
Let us take $p(x) = B_n(x) \in \mathbb{P}_n$, $(n \geq 0)$. Then, by Theorem 8, we get
\begin{equation}\begin{split}\label{67}
B_n(x) = \sum_{m=0}^n a_m b_{m,\lambda }(x),\,\,(n \geq 0),
\end{split}\end{equation}
where
\begin{equation}\begin{split}\label{68}
a_m &= \frac{1}{m!} < (2-e^t)(\tfrac{1}{\lambda }(e^{\lambda t}-1))^m | B_n(x) >\\
&=\frac{1}{m!}<(2-e^t)\lambda^{-m}m! \sum_{k=m}^\infty S_2(k,m) \lambda ^k \frac{t^k}{k!} | B_n(x)>\\
&= \sum_{k=m}^n \lambda ^{k-m}S_2(k,m) \frac{1}{k!} <2-e^t | t^k B_n(x) > \\
&= \sum_{k=m}^n \lambda ^{k-m}S_2(k,m) {n \choose k} <2-e^t| B_{n-k}(x) > \\
&=\sum_{k=m}^n \lambda ^{k-m}S_2(k,m) {n \choose k}\big( 2B_{n-k}-B_{n-k}(1) \big)\\
&=\sum_{k=m}^n \lambda ^{k-m}S_2(k,m) {n \choose k}\big(B_{n-k}- \delta_{1,n-k} \big)\\
&= \sum_{k=m}^n \lambda ^{k-m} S_2(k,m) {n \choose k} B_{n-k} - n \lambda ^{n-m-1} S_2(n-1,m).
\end{split}\end{equation}
Therefore, by \eqref{67} and \eqref{68}, we obtain the following theorem.
\begin{thm}
For $n \geq 0$, we have
\begin{equation*}\begin{split}
B_n(x) &= \sum_{m=0}^n \sum_{k=m}^n \lambda ^{k-m} S_2(k,m) {n \choose k} B_{n-k} b_{m,\lambda }(x) \\
&\,\, - n \sum_{m=0}^{n-1} \lambda ^{n-m-1} S_2(n-1,m) b_{m,\lambda }(x).
\end{split}\end{equation*}
\end{thm}
For $b_{n,\lambda }(x) \sim ( 2-e^t, \tfrac{1}{\lambda }(e^{\lambda t}-1))$, $(x)_{n,\lambda } \sim (1, \tfrac{1}{\lambda }(e^{\lambda t}-1))$, $(n \geq 0)$, we have
\begin{equation}\begin{split}\label{69}
(x)_{n,\lambda } = \sum_{m=0}^n C_{n,m}b_{m,\lambda }(x),\,\,(n \geq 0),
\end{split}\end{equation}
where
\begin{equation}\begin{split}\label{70}
C_{n,m}&= \frac{1}{m!}< (2-e^{\frac{1}{\lambda }\log(1+\lambda t)})^{-1}(\tfrac{1}{\lambda }(e^{\log(1+\lambda t)}-1)^m | x^n >\\
&=\frac{1}{m!} < (2-(1+\lambda t)^{\frac{1}{\lambda }})^{-1}t^m |x^n>\\
&= {n \choose m} <\sum_{k=0}^{\infty} b_{k,\lambda} \frac{t^k}{k!}|x^{n-m}>\\
&= {n \choose m}b_{n-m,\lambda}.
\end{split}\end{equation}
Therefore, by \eqref{70}, we obtain the following theorem.
\begin{thm}
For $n \geq 0$, we have
\begin{equation*}\begin{split}
(x)_{n,\lambda } = \sum_{m=0}^n {n \choose m}b_{n-m,\lambda}b_{m,\lambda}(x).
\end{split}\end{equation*}
\end{thm}

The Korobov polynomials $K_{n,\lambda }(x)$ for $\lambda  \neq 0,1$ are also defined in terms of generating function by
\begin{equation}\begin{split}\label{71}
\frac{\lambda t}{(1+t)^\lambda -1}(1+t)^x = \sum_{n=0}^\infty K_{n,\lambda }(x) \frac{t^n}{n!},\quad (\textnormal{see} \,\, [6]).
\end{split}\end{equation}
Thus, by \eqref{71}, we get
\begin{equation}\begin{split}\label{72}
K_{n,\lambda }(x) \sim \left( \frac{e^{\lambda t}-1}{\lambda (e^t-1)}, e^t-1\right),
\end{split}\end{equation}
For $K_{n,\lambda }(x) \sim \left( \frac{e^{\lambda t}-1}{\lambda (e^t-1)}, e^t-1 \right)$, $b_{n,\lambda }(x) \sim \left( 2-e^t, \frac{1}{\lambda }(e^{\lambda t}-1)\right)$,
we have
\begin{equation}\begin{split}\label{73}
K_{n,\lambda }(x)= \sum_{m=0}^n C_{n,m}b_{m,\lambda }(x),\,\,(n \geq 0),
\end{split}\end{equation}
where
\begin{equation}\begin{split}\label{74}
C_{n,m}&= \frac{1}{m!} \Big< (1-t)\frac{\lambda t}{(1+t)^{\lambda}-1}\Big(\frac{1}{\lambda}\big((1+t)^{\lambda}-1\big)\Big)^m \Big| x^n \Big>\\
&=\frac{\lambda^{-m}}{m!} \sum_{l=0}^m {m \choose l} (-1)^{m-l} \Big< (1-t)\frac{\lambda t}{(1+t)^{\lambda}-1}(1+t)^{l\lambda} \Big| x^n\Big>\\
&=\frac{\lambda^{-m}}{m!} \sum_{l=0}^m {m \choose l} (-1)^{m-l} \Big< (1-t)  \Big|\sum_{j=0}^{n}K_{j,\lambda}(l\lambda)\frac{t^j}{j!} x^n\Big>\\
&=\frac{\lambda^{-m}}{m!} \sum_{l=0}^m {m \choose l} (-1)^{m-l}\sum_{j=0}^{n}K_{j,\lambda}(l\lambda){n \choose j}(\delta_{j,n}-\delta_{j,n-1})\\
&=\frac{\lambda^{-m}}{m!} \sum_{l=0}^m {m \choose l} (-1)^{m-l}\big(K_{n,\lambda}(l \lambda)-nK_{n-1,\lambda}(l\lambda)\big).
\end{split}\end{equation}
Therefore, by \eqref{73} and \eqref{74}, we obtain the following theorem.
\begin{thm}
For $n \geq 0$, we have
\begin{equation*}\begin{split}
K_{n,\lambda }(x) = \sum_{m=0}^n \left( \frac{\lambda^{-m}}{m!} \sum_{l=0}^m {m \choose l} (-1)^{m-l}\big(K_{n,\lambda}(l \lambda)-nK_{n-1,\lambda}(l\lambda)\big) \right) b_{m,\lambda }(x).
\end{split}\end{equation*}
\end{thm}

\end{document}